\documentclass{article}

\usepackage{arxiv}
\usepackage{graphicx}
\usepackage{bbold}
\usepackage{optidef}
\usepackage{multirow}

\usepackage{amsthm}
\theoremstyle{definition}
\newtheorem{defn}{Definition}
\theoremstyle{plain}
\newtheorem{thm}{Theorem}

\newtheorem{prop}[thm]{Proposition}

\usepackage[utf8]{inputenc}
\usepackage[T1]{fontenc}
\usepackage{optidef}
\usepackage{hyperref}
\usepackage{nicematrix}
\usepackage{xfrac} 
\usepackage{caption}
\usepackage{subcaption}
\usepackage{comment}
\usepackage[ruled,linesnumbered]{algorithm2e}

\SetCommentSty{mycommfont}
\hypersetup{
    colorlinks=true,
    linkcolor=blue,
    filecolor=magenta,      
    urlcolor=black,
}
\usepackage[numbers]{natbib}

\usepackage[utf8]{inputenc} 
\usepackage[T1]{fontenc}    
\usepackage{hyperref}       
\usepackage{url}            
\usepackage{booktabs}       
\usepackage{amsmath}
\usepackage{amsfonts}       
\usepackage{nicefrac}       
\usepackage{microtype}      
\usepackage{lipsum}		
\usepackage{graphicx}
\usepackage{doi}

\newcommand{\argmax}{\mathop{\rm arg~max}\limits}
\newcommand{\argmin}{\mathop{\rm arg~min}\limits}
\title{A Particle-Based Algorithm for Distributional Optimization on \textit{Constrained Domains} via Variational Transport and Mirror Descent}

\author{ \href{https://orcid.org/0000-0003-0380-4197}{\includegraphics[scale=0.06]{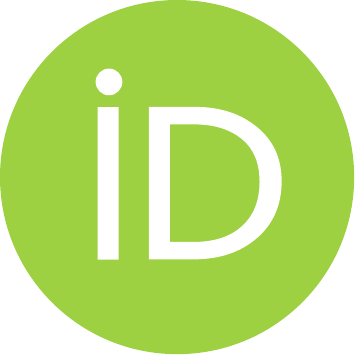}\hspace{1mm}Dai Hai Nguyen}\\
	Department of Computer Science\\
	University of Tsukuba, Japan\\
	\texttt{hai@cs.tsukuba.ac.jp} \\
    	\And
	\href{https://orcid.org/0000-0002-5789-7547}{\includegraphics[scale=0.06]{orcid.pdf}\hspace{1mm}Tetsuya Sakurai} \\
    Department of Computer Science\\
	University of Tsukuba, Japan\\
	\texttt{sakurai@cs.tsukuba.ac.jp} \\
}

\renewcommand{\shorttitle}{\textbf{mirrorVT}: A Distributional Optimization Framework on Constrained Domains}

\hypersetup{
pdftitle={A template for the arxiv style},
pdfsubject={q-bio.NC, q-bio.QM},
pdfauthor={David S.~Hippocampus, Elias D.~Striatum},
pdfkeywords={First keyword, Second keyword, More},
}

\begin{document}
\maketitle

\begin{abstract}
We consider the optimization problem of minimizing an objective functional, which admits a variational form and is defined over probability distributions on the \textit{constrained domain}, which poses challenges to both theoretical analysis and algorithmic design. Inspired by the mirror descent algorithm for constrained optimization, we propose an iterative particle-based algorithm, named Mirrored Variational Transport (\textbf{mirrorVT}), extended from the Variational Transport framework \cite{liu2021infinite} for dealing with the constrained domain. In particular, for each iteration, \textbf{mirrorVT} maps particles to an unconstrained dual domain induced by a mirror map and then approximately perform Wasserstein gradient descent on the manifold of distributions defined over the dual space by pushing particles. At the end of iteration, particles are mapped back to the original constrained domain. Through simulated experiments, we demonstrate the effectiveness of \textbf{mirrorVT} for minimizing the functionals over probability distributions on the simplex- and Euclidean ball-constrained domains. We also analyze its theoretical properties and characterize its convergence to the global minimum of the objective functional.
\end{abstract}

\section{Introduction}
Many problems in machine learning and computational statistics involve a distributional optimization problem in which we wish to optimize a functional $F:\mathcal{P}_{2}(\mathcal{X})\rightarrow \mathbb{R}$ of probability distributions : $\min_{p\in \mathcal{P}_{2}(\mathcal{X})}F(p)$ where $\mathcal{P}_{2}(\mathcal{X})$ denote the set of probability distributions defined over the domain $\mathcal{X}$ with finite second-order moments. A common class of methods falling into such a category is gradient-based Markov Chain Monte Carlo (MCMC) sampling, commonly used in Bayesian inference (see \cite{welling2011bayesian, wibisono2018sampling, ma2015complete} and references therein).
These sampling methods attempt to approximate a target distribution $p^{*}$ by generating a set of particles and updating the particles in such a way that a dissimilarity $D$ between the approximate empirical probability measure $p$ given by the current set of particles and the target $p^{*}$ is minimized. The dissimilarity function $D(p,p^{*})$ includes Kullback-Leiber (KL) divergence, Jensen-Shanon (JS) divergence or Wasserstein distance in the optimal transport \cite{villani2009optimal}. By setting $F(p)=D(p,p^{*})$, the sampling task can be seen as a distributional optimization problem.

In general, this optimization problem can be solved by an iterative algorithm, named \textit{Wasserstein Gradient Descent} \cite{zhang2016first}, which has two steps for each iteration: 1) obtain Wasserstein gradient of $F$ with respect to the current probability distribution and 2) perform the exponential mapping on $\mathcal{P}_{2}(\mathcal{X})$. Variational Transport (\textbf{VT}, \cite{liu2021infinite}) algorithm is recently proposed to make it useful for the distributional optimization by approximating a probability distribution using a set of particles and updating the particles to minimize $F$. To this end, \textbf{VT} assumes that $F$ admits a variational form and solve a \textit{variational maximization problem} to approximate Wasserstein gradient of $F$ and then use the obtained solution to specify a direction to update each particle. This can be seen as a forward discretization of the Wasserstein gradient flow \cite{santambrogio2017euclidean}.

However, applying \textbf{VT} to constrained distributions remains a challenge as it can push the particles outside the domain by following the direction given by the solution of the variational maximization problem. In this paper, we address this challenge by proposing mirror Variational Transport (\textbf{mirrorVT}) working on the dual space with unconstrained domain instead of the (original) primal space with constrained domain via a mirror map. Our approach is motivated by the Mirror Descent Algorithm (\textbf{MD}, \cite{beck2003mirror}) for the constrained optimization problem.

We summarize our contributions of this paper as follows. 1) we extend the variational transport framework to deal with the distributional optimization problem on the \textit{constrained domain}. Motivated by the \textbf{MD} algorithm, we convert the original optimization problem from the primal space (constrained domain) into dual space (unconstrained domain) via the mirror map. 2) We analyze theoretical properties and the convergence of \textbf{mirrorVT} to the optimal solution of the objective functional. 3) We conduct simulated experiments on \textit{simplex- and Euclidean ball-constrained domains} to demonstrate the effectiveness of \textbf{mirrorVT} in dealing with the distributional optimization on the constrained domain.

\section{Related Works}
Our work is relevant to the line of research on sampling methods, particularly gradient-based MCMC, which is
a class of particle-based sampling algorithms for a target distribution. For example, see \cite{wibisono2018sampling,welling2011bayesian,xu2018global} and references therein. 
Our work is more related to \cite{cheng2018convergence, ma2015complete} in which the finite-time convergence of gradient-based methods in terms of KL-divergences in $\mathcal{P}_{2}(\mathcal{X})$ is comprehensively studied.

In addition to gradient-based MCMC, our method is also closely related to \textbf{VT} \cite{liu2021infinite}, which utilizes the optimal transport framework and the variational form of the objective functional for solving the distributional optimization problems via particle approximation. In particular, in each iteration, \textbf{VT} estimates the Wasserstein gradient by solving a variational maximization problem associated with the objective functional, and then performs Wasserstein gradient descent by simply pushing the particles. By utilizing the variational representation of the objective functional, \textbf{VT} can be applied to optimize functionals beyond the KL-divergences commonly targeted by gradient-based MCMC algorithms. However, when working on the constrained domain, it can push the particles \textit{outside the domain}, making it infeasible to preserve the constrained domain.

Another line of related research in dealing with sampling for constrained domain is inspired by the classical \textbf{MD} \cite{beck2003mirror}, for example, see \cite{ahn2021efficient,hsieh2018mirrored} and references therein. The basic idea of these works is to transform the constrained sampling problems into unconstrained ones by using a \textit{mirror map}. The finite-time convergence of these methods is also well studied.
Furthermore, our work is most related to recently proposed Stein Variational Mirror Descent (\textbf{SVMD}, \cite{shi2021sampling}), which is a variant of Stein Variational Gradient Desccent (\textbf{SVGD}, \cite{liu2016stein}) and suitable for sampling from the constrained domain and non-Euclidean geometries. It basically minimizes the KL-divergence to the constrained target distribution by pushing particles in a dual space induced by a mirror map. In this paper, we exploit the connection between our method and \textbf{SVMD}, and show that our method relates to \textbf{SVMD} when we utilize the KL-divergence as the objective functional and approximate the optimal direction to push the particles via the integral operator \cite{rosasco2009note}. Furthermore, while \textbf{SVMD} can be seen as a distributional optimization problem with the KL-divergence as the objective functional, our method can be applied to a broad class of functionals due to the use of the \textbf{VT} framework.

\section{Background: Variational Transport and Mirror Descent}
We consider the following distributional optimization problem on a constrained domain:
\begin{align}
\label{eqn:distopt}
    \min_{p\in \mathcal{P}_{2}(\mathcal{X})}F(p)
\end{align}
where $\mathcal{X}$ denotes a $d$-dimensional Riemannian manifold with the constrained support. In this section, we briefly describe preliminaries on optimal transport, the \textbf{VT} and \textbf{MD} algorithms which are closely related to our proposed algorithm.

\subsection{Preliminaries on Optimal Transport and Wasserstein Space}
Given a measurable map $T:\mathcal{X}\rightarrow \mathcal{X}$ and $p\in \mathcal{P}_{2}(\mathcal{X})$, we say that $q$ is the \textit{push-forward measure} of $p$ under $T$, denoted by $q=T\sharp p$, if for every Borel set $E\subseteq \mathcal{X}$, $q(E)=p(T^{-1}(E))$. For any $p,q\in \mathcal{P}_{2}(\mathcal{X})$, the $2$-Wasserstein distance $\mathcal{W}_{2}(p,q)$ is defined as:
\begin{align}
   \mathcal{W}_{2}^{2}(p, q)= \inf_{\pi \in \Pi(p,q)} \int_{\mathcal{X}\times \mathcal{X}} \lVert \textbf{x} - \textbf{x}^\prime\rVert_{2}^{2}\mathrm{d}\pi(\textbf{x},\textbf{x}^\prime)
\end{align}
where $\Pi(p,q)$ is all probability measures on $\mathcal{X}\times \mathcal{X}$ whose two marginals are equal to $p$ and $q$, $\lVert\cdot \rVert_{2}$ denotes the Euclidean norm. It is known that the metric space $(\mathcal{P}_{2}(\mathcal{X}), \mathcal{W}_{2})$, also known as Wasserstein space, is an infinite-dimensional geodesic space \cite{villani2009optimal}.

Given a functional $F:\mathcal{P}_{2}(\mathcal{X})\rightarrow \mathbb{R}$, the first variation of $F$ evaluated at $p$, denoted by $\partial F(p)/\partial p:\mathcal{X}\rightarrow \mathbb{R}$, is  given as follows:
\begin{align}
    \lim_{\epsilon\rightarrow 0}\frac{1}{\epsilon} \left( F(p+\epsilon \chi) - F(p) \right)=\int_{\mathcal{X}} \frac{\partial F(p)}{\partial p}(\textbf{x})\chi(\textbf{x})\mathrm{d}(\textbf{x})
\end{align}
for all $\chi=q-p$, where $q\in \mathcal{P}_{2}(\mathcal{X})$. With mild regularity assumptions, the Wasserstein gradient of $F$, denoted by $\texttt{grad}F$, relates to the gradient of the first variation of $F$ via the following continuity equation:
\begin{align}
    \texttt{grad}F(p)(\textbf{x})=-\texttt{div}\left(p(\textbf{x})\nabla\frac{\partial F(p)}{\partial p}(\textbf{x})\right)
\end{align}
for all $\textbf{x}\in \mathcal{X}$, where $\texttt{div}$ denotes the divergence operator. For every vector field $v:\mathcal{X}\rightarrow\mathcal{X}$, we have:
\begin{align}
\label{eqn:div}
    \int_{\mathcal{X}}\langle \nabla\frac{\partial F(p)}{\partial p}(\textbf{x}),v(\textbf{x}) \rangle p(\textbf{x})\mathrm{d}\textbf{x}=- \int_{\mathcal{X}}\frac{\partial F(p)}{\partial p}(\textbf{x})\texttt{div}(p(\textbf{x})v(\textbf{x}))\mathrm{d}\textbf{x}
\end{align}
We refer the readers to \cite{santambrogio2017euclidean} for more details.

\subsection{The Variational Transport (VP) Algorithm}
When $\mathcal{X}$ is the unconstrained domain, the functional gradient descent with respect to the geodesic distance can be used to directly optimize $F$. Essentially, it constructs a sequence of probability distributions $\left\{ p_{t}\right\}_{t\geq 1}$ in $\mathcal{P}_{2}(\mathcal{X})$ as follows:
\begin{align}
\label{eqn:functionalgd}
    p_{t+1}\leftarrow \texttt{Exp}_{p_{t}}\left\{ -\eta_{t}\cdot \texttt{grad}F(p_{t})\right\}
\end{align}
where $\eta_{t}$ is the stepsize, $\texttt{grad}F(p_{t})$ is the Wasserstein gradient evaluated at $p_{t}$ and $\texttt{Exp}_{p}$ denotes the exponential mapping by which the probability measure $p$ is moved along a given direction on $\mathcal{P}_{2}(\mathcal{X})$. See \cite{zhang2016first} for more details.

The Variational Transport (\textbf{VP}, \cite{liu2021infinite}) algorithm is introduced for solving the distributional optimization problem (\ref{eqn:distopt}) by approximating $p_{t}$ by an empirical measure $\Tilde{p}_{t}$ of $N$ particles $\left\{\textbf{x} _{t,i} \right\}_{i\in \left[N\right]}$ and assuming that the functional $F$ admits the following variational form:
\begin{align}
\label{eqn:varform}
    F(p_{t}) = \sup_{f\in\mathcal{F}}{\left\{\int{f(\textbf{x})p_{t}(\textbf{x})\mathrm{d}\textbf{x}}-F^{*}(f)\right\}}
\end{align}
where $\mathcal{F}$ is a class of square-integrable functions on $\mathcal{X}$ with respect to the Lebesgue measure and $F^{*}:\mathcal{F}\rightarrow \mathbb{R}$ is a convex functional of $F$. The advantage of the variational functional objective is that the Wasserstein gradient can be calculated from the solution $f^{*}_{t}$ to the problem (\ref{eqn:varform}), which can be estimated using samples from $p_{t}$. Specifically, it is shown that $f^{*}_{t}=\partial F/\partial p_{t}$, which is the first variation of $F$ (see Proposition 3.1 in \cite{liu2021infinite}). Furthermore,  under the assumption that $\nabla f^{*}_{t}$ is $h$-Lipschitz continuous, then for any $\eta_{t}\in \left[ 0,1/h\right)$, the exponential mapping in (\ref{eqn:functionalgd}) is shown to be equivalent to the push-forward mapping defined by $f^{*}_{t}$. That is for $\textbf{x}_{t,i}$ drawn from $p_{t}$:

\begin{align}
    \textbf{x}_{t+1,i}\leftarrow \texttt{Exp}_{\textbf{x}_{t,i}}\left\{ -\eta_{t}\cdot \nabla f^{*}_{t}(\textbf{x}_{t,i})\right\}
\end{align}
where $\textbf{x}_{t+1,i}$ is the updated particle which is drawn from $p_{t+1}$, $\texttt{Exp}_{\textbf{x}}(\eta\cdot \nabla u)$ denotes the transportation map which sends $\textbf{x}\in  \mathcal{X}$ to a point $\textbf{x}+\eta\cdot \nabla u \in \mathcal{X}$ (see Proposition 3.2 in \cite{liu2021infinite}).

In addition, \textbf{VP} estimates the solution $f^{*}_{t}$ by solving the empirical variational maximization problem with finite samples drawn from $p_{t}$ via stochastic gradient descent on the domain $\mathcal{X}$ (see Algorithm \ref{alg:vfm}):

\begin{align}
\label{eqn:varformmax}
     \Tilde{f_{t}^{*}}= \argmax_{f\in\mathcal{F}}{\left\{\frac{1}{N}\sum_{i=1}^{N}f(\textbf{x}_{t,i})-F^{*}(f)\right\}}
\end{align}

where $\mathcal{F}$ is a function class, which can be specified to be the following class of deep neural networks:
\begin{align}
    \label{eqn:nnclass}
    \Tilde{\mathcal{F}}=\left\{\Tilde{f}| \Tilde{f}(\textbf{x})=\frac{1}{\sqrt{n_{w}}}\sum_{i=1}^{n_{w}}b_{i}\sigma([\textbf{w}]_{i}^{T}\textbf{x})\right\}
\end{align}
where $n_{w}$ is the width of the neural networks, $[\textbf{w}]_{i}\in \mathbb{R}^{d}$, $\textbf{w}=([\textbf{w}]_{1},...,[\textbf{w}]_{n_{w}})^{T}\in \mathbb{R}^{n_{w}\times d}$ is the input weight, $\sigma$ denotes a smooth activation function, and $b_{i}\in \left\{ -1,1\right\}$. In addition, in each iteration, the weights $\textbf{w}$ is guaranteed to lie in the $l_{2}$-ball centered at the initial weights $\textbf{w}(0)$ with radius $r_{f}$ defined as $\mathcal{B}^{0}(r_{f})=\left\{\textbf{w}: \lVert\textbf{w} - \textbf{w}(0) \rVert_{2}\leq r_{f} \right\}$. We refer the readers to \cite{liu2021infinite} for more details of this neural network parameterization. This choice of neural network class facilitates the analysis of the gradient error induced by the difference between $f^{*}_{t}$ and $ \Tilde{f_{t}^{*}}$.

However, when the domain $\mathcal{X}$ is constrained, e.g. $d$-dimensional simplex, the updates of \textbf{VP} may fail to preserve the constrained domain by pushing particles outside the domain. We tackle this problem in this paper by adopting the principle of the following mirror descent algorithm.

\begin{table}[]
    \centering
    \caption{Common choice  of functionals $F$, $G$ and their corresponding conjugate functionals}
     \resizebox{\textwidth}{!}  
     {
    \begin{tabular}{l l l l l l}
    \hline
    Divergence & Functional $F(p)$ & Functional $G(q)$ & Conjugate $F^{*}(f)$ & Conjugate $G^{*}(g)$\\
    \hline
    KL  & $\texttt{KL}(p(\textbf{x})||p^{*}(\textbf{x}))$ & $\texttt{KL}(q(\textbf{y})||q^{*}(\textbf{y}))$ & $\texttt{log}\mathbb{E}_{p^{*}} e^{f(\textbf{x})}$ & $\texttt{log} \mathbb{E}_{q^{*}}e^{g(\textbf{y})}$\\ 
    JS  & $\text{JS}((p(\textbf{x})||p^{*}(\textbf{x}))$ & $\text{JS}(q(\textbf{y})||q^{*}(\textbf{y}))$ & $-\frac{1}{2}E_{p^{*}} \texttt{log}\left(1-2e^{2f(\textbf{y})}\right)-\frac{1}{2}\texttt{log}2$ & $-\frac{1}{2}E_{q^{*}} \texttt{log}\left(1-2e^{2g(\textbf{y})}\right)-\frac{1}{2}\texttt{log}2$\\ 
    1-Wasserstein & $\mathcal{W}_{1}(p(\textbf{x})||p^{*}(\textbf{x}))$ & $\frac{1}{\alpha}\mathcal{W}_{1}(q(\textbf{\textbf{y}})||q^{*}(\textbf{y}))$ & $\mathbb{E}_{p^{*}}\left[f(\textbf{x})+\{\lVert f \rVert_{L}\leq 1 \} \right]$ & $\mathbb{E}_{q^{*}}\left[g(\textbf{y})+\{\lVert g \rVert_{L}\leq \frac{1}{\alpha} \} \right]$\\ 
    \end{tabular}
    }
    \label{tab:examples}
\end{table}

\subsection{The Mirror Descent (MD) Algorithm}
Gradient Descent is the standard algorithm to optimize a target function $f$ on the unconstrained domain $\mathcal{X}$ by solving the following optimization problem for each step $t$:
\begin{align}
\label{eqn:gd}
     \textbf{x}_{t+1}=\argmin_{\textbf{x}\in\mathcal{X}}{\langle\nabla f(\textbf{x}_{t}),\textbf{x}\rangle + \frac{1}{2\eta_{t}}\lVert \textbf{x}-\textbf{x}_{t}\rVert^{2}_{2}}
\end{align}
To deal with the constrained optimization problems, the Mirror Descent (\textbf{MD}) algorithm replaces $\lVert \cdot\rVert_{2}$ in (\ref{eqn:gd}) with a function $\varphi$ that reflects the geometry of the problem \cite{beck2003mirror}. The \textbf{MD} algorithm chooses $\Phi$ to be the Bregman divergence induced by a strongly convex function $\varphi:\mathcal{X}\rightarrow \mathbb{R}$ as follows: $\Phi(\textbf{x}^\prime, \textbf{x})=\varphi(\textbf{x}^\prime) - \varphi(\textbf{x}) - \langle \nabla \varphi(\textbf{x}),\textbf{x}^\prime-\textbf{x}\rangle$ for $\textbf{x}^\prime, \textbf{x}\in \mathcal{X}$. Then, the solution of (\ref{eqn:gd}) for each step becomes:
\begin{align}
\label{eqn:mdupdate}
     \textbf{x}_{t+1}=\nabla \varphi^{*} \left( \nabla \varphi(\textbf{x}_{t})-\eta_{t} \nabla f(\textbf{x}_{t}) \right)
\end{align}
where $\varphi^{*}(\textbf{y})=\sup_{\textbf{x} \in \mathcal{X}}{\langle \textbf{x},\textbf{y}\rangle -\varphi(\textbf{x})}$ is the convex conjugate of function $\varphi$ and $\nabla \varphi^{*}(\textbf{y})=(\nabla \varphi)^{-1}(\textbf{y})$ is the inverse map. Intuitively, the \textbf{MD} update (\ref{eqn:mdupdate}) is composed of three steps: 1) mapping $\textbf{x}_{t}$ to $\textbf{y}_{t}$ by $\nabla \varphi$, 2) applying the update: $\textbf{y}_{k+1}=\textbf{y}_{t} - \eta_{t}\nabla f(\textbf{x}_{t})$, and 3) mapping back through $\textbf{x}_{t+1} = \nabla \varphi^{*}(\textbf{y}_{t+1})$.

\section{The Mirrored Variational Transport (mirrorVT) Algorithm}
In what follows, we introduce the main algorithm to tackle the distributional optimization problem (\ref{eqn:distopt}) when the domain $\mathcal{X}$ is constrained.
\subsection{From The Primal Space To Dual Space}
Motivated by the \textbf{MD} algorithm for constrained optimization problem, we propose to solve the problem (\ref{eqn:distopt}) in the dual space with unconstrained support via a $\alpha$-strongly convex function $\varphi$. In particular, the variational form of $F$ in (\ref{eqn:varform}) can be written as follows:

\begin{align}
\begin{split}
    \label{eqn:FandG}
    F(p) &= \sup_{f\in\mathcal{F}}{\left\{\int{f(\textbf{x})p(\textbf{x})\mathrm{d}\textbf{x}}-F^{*}(f)\right\}}\\
         &= \sup_{f\in\mathcal{F}}{\left\{\int{f(\nabla \varphi^{*}(\textbf{y}))q(\textbf{y})\mathrm{d}y}-F^{*}(f\circ \nabla \varphi^{*})\right\}}\\
         &= \sup_{g\in\mathcal{G}}{\left\{\int{g(\textbf{y})q(\textbf{y})\mathrm{d}\textbf{y}}-G^{*}(g)\right\}}=G(q)\\
\end{split}
\end{align}
where $\textbf{y}=\nabla \varphi(\textbf{x}) \in \mathcal{Y}$, $\mathcal{Y}$ denotes the $d$-dimensional Riemannian manifold with the unconstrained support, $q=(\nabla \varphi)\sharp p$ is the push-forward measure of $p$ induced by the mapping $\nabla \varphi$, 
$g=f\circ\nabla \varphi^{*} \in \mathcal{G}$, 
$\mathcal{G}$ denotes a function class defined on the domain $\mathcal{Y}$. 
To be more intuitive, several examples of $F^{*}$, $G^{*}$ are shown in Table \ref{tab:examples}.

\begin{figure}[t]
	\centerline{\includegraphics[width=1.0\columnwidth]{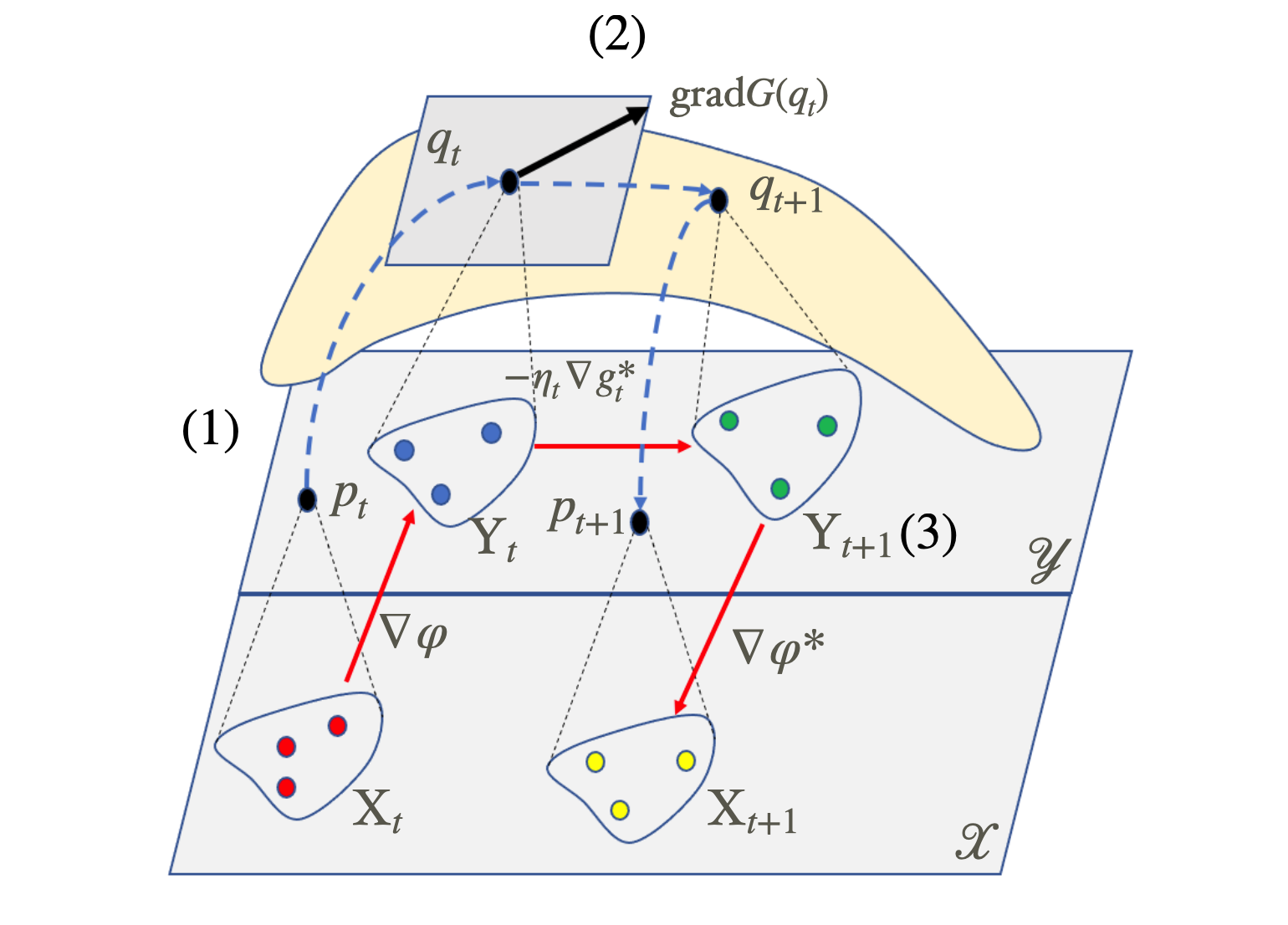}}
	\caption{Overview of Mirror Variational Transport (\textbf{mirrorVT}). In $t$-th iteration, there are three main steps: (1) \textbf{mirrorVT} maps the particle set $\textbf{X}_{t}$ to $\textbf{Y}_{t}$ in the unconstrained domain $\mathcal{Y}$ via the mirror map $\nabla \varphi$, (2) \textbf{mirrorVT} pushes the particle set $\textbf{Y}_{t}$ along the direction $\nabla g^{*}_{t}$ in the domain $\mathcal{Y}$ to obtain the updated particle set $\textbf{Y}_{t+1}$, and (3) \textbf{mirrorVT} maps $\textbf{Y}_{t+1}$ back to $\textbf{X}_{t+1}$ in the constrained domain $\mathcal{X}$ via the map $\nabla \varphi^{*}$.}
	\label{fig:overview}
\end{figure}

The advantages of defining the functional $G$ on the unconstrained domain $\mathcal{Y}$ are as follows: 1) $G$ is defined over the unconstrained domain $\mathcal{Y}$, which enables to apply \textbf{VP} to optimize $G$ in the dual space; 2) It is straightforward to verify the following connection between $F$ and $G$: $\partial G/\partial q=g_{q}^{*}$ and $\nabla g_{q}^{*}(\textbf{y})= \nabla^{2}\varphi(\textbf{x})^{-1}\nabla f_{p}^{*}(\textbf{x})$, where $f_{p}^{*}$ and $g_{q}^{*}$ are the optimal solutions to the variational maximization problem (\ref{eqn:FandG}) of $F$ and $G$, respectively. 

So we propose the mirrored variational transport (\textbf{mirrorVT}) algorithm, as illustrated in Fig. \ref{fig:overview}, to optimize the functional $F$. Initially, a set of particles $\textbf{X}_{1}=\left\{\textbf{x}_{1,i} \right\}_{i\in [N]}$ is obtained by drawing $N$ i.i.d observations from the initial probability distribution $p_{1}$. In each iteration $t$, we maintain two set of $N$ particles $\textbf{X}_{t}=\left\{\textbf{x}_{t,i} \right\}_{i\in [N]}$ and $\textbf{Y}_{t}=\left\{\textbf{y}_{t,i} \right\}_{i\in [N]}$ as follows. First, \textbf{mirrorVT} compute the solution $f^{*}_{t}$ to (\ref{eqn:varform}) 
based on the current particle set $\textbf{X}_{t}$. Then, \textbf{mirrorVT} maps the current particle set $\textbf{X}_{t}$ to the domain $\mathcal{Y}$ via the mirror map $\nabla \varphi$ (see step 1 of Fig. \ref{fig:overview}):
\begin{align}
    \textbf{y}_{t,i} \leftarrow \nabla \varphi(\textbf{x}_{t,i})
\end{align}
for all $i\in [N]\footnote{We use [N] to indicate the list $\left[1,2,...,N\right]$ throughout the rest of the paper}$. This is to transform the optimization problem from constrained domain into the unconstrained one ($\mathcal{Y}$) and equivalent to pushing the empirical measure $p_{t}=1/N\sum_{i\in[N]}\delta_{\textbf{x}_{t,i}}$ by the push-forward measure induced by $\nabla \varphi$: $q_{t}=\left[\nabla \varphi \right]\sharp p_{t}$. 

With the obtained $f^{*}_{t}$, we assume that $\nabla^{2}\varphi(\textbf{x})^{-1}\nabla f^{*}_{t}(\textbf{x})$ is $h$-Lipschitz continuous and one may verify that $\nabla g^{*}_{t}(\textbf{y})$ is $h/\alpha$-Lipschitz continuous\footnote{For any $\textbf{x},\textbf{x}^\prime \in \mathcal{X}$ and $\textbf{y}=\nabla \varphi(\textbf{x}),\textbf{y}^\prime=\nabla \varphi(\textbf{x}^\prime)$, we have: $\lVert \nabla g^{*}_{t}(\textbf{y})-\nabla g^{*}_{t}(\textbf{y}^\prime) \rVert_{2} = \lVert \nabla^{2} \varphi(\textbf{x})^{-1}\nabla f^{*}_{t}(\textbf{x})-\nabla^{2} \varphi(\textbf{x}^\prime)^{-1}\nabla f^{*}_{t}(\textbf{x}^\prime)\rVert_{2} \leq h \lVert \textbf{x}-\textbf{x}^{\prime} \rVert_{2} = h \lVert \nabla \varphi^{*}(\textbf{y})-\nabla \varphi^{*}(\textbf{y}^\prime) \rVert_{2} \leq h/\alpha\lVert \textbf{y}-\textbf{y}^\prime \rVert_{2}$, where the last inequality holds as $\varphi^{*}$ is $1/\alpha$-smooth.}. Hence for any $\eta_{t}\in \left[0, \alpha/h\right)$, we can enable \textbf{mirrorVT} to push the particles in the dual domain as follows (see step 2 of Fig. \ref{fig:overview}):
\begin{align}
    \textbf{y}_{t+1,i}\leftarrow \texttt{Exp}_{\textbf{y}_{t,i}}\left\{ -\eta_{t} \nabla ^{2}\varphi(\textbf{x}_{t,i})^{-1}f^{*}_{t}(\textbf{x}_{t,i})\right\}
\end{align}
for all $i\in [N]$. This is equivalent to updating the empirical measure $q_{t}=1/N\sum_{i\in[N]}\delta_{\textbf{y}_{t,i}}$ by the push-forward measure $\left[\texttt{Exp}_{\mathcal{Y}}(-\eta_{t} \nabla g^{*}_{t}) \right]\sharp q_{t}$, which is realized by applying Proposition 3.1 in \cite{liu2021infinite} for the function $g^{*}_{t}$ defined over the dual domain $\mathcal{Y}$. At the end of $t$-th iteration, \textbf{mirrorVT} maps particles back to the original constrained domain $\mathcal{X}$ (see step 3 of Fig. \ref{fig:overview}):
\begin{align}
    \textbf{x}_{t+1,i} \leftarrow \nabla \varphi^{*}(\textbf{y}_{t+1,i})
\end{align}
for all $i\in [N]$. This step is equivalent to pushing the updated empirical measure $q_{t+1}=1/N\sum_{i\in[N]}\delta_{\textbf{y}_{t+1,i}}$
from the dual domain to the primal one by the push-forward measure $\left[\nabla \varphi^{*} \right]\sharp q_{t+1}$.

\subsection{Continuous-time dynamics of $p_{t}$}
We can view steps of \textbf{mirrorVT} as a discretization of the following continuous-time dynamics of particles (as $\eta_{t}\rightarrow 0$):
\begin{equation}
    \mathrm{d}\textbf{y}_{t}=-\nabla^{2}\varphi(\textbf{x}_{t})^{-1}\nabla f^{*}_{t}(\textbf{x}_{t})\mathrm{d}t \text{ and } \textbf{x}_{t}=\nabla \varphi^{*}(\textbf{y}_{t})\mathrm{d}t
\end{equation}
which is equivalent to
\begin{equation}
    \label{eqn:dxdt}
    \frac{\mathrm{d}\textbf{x}_{t}}{\mathrm{d}t}=\left( \nabla_{\textbf{x}_{t}}\textbf{y}_{t}\right)^{-1}\frac{\mathrm{d}\textbf{y}_{t}}{\mathrm{d}t}=- \nabla^{2}\varphi(\textbf{x}_{t})^{-2}\nabla f^{*}_{t}(\textbf{x}_{t})
\end{equation}
by applying the chain rule. We now consider the continuous-time dynamic of $p_{t}$ and its limit. Let $\gamma_{t}:\mathbb{R}^{+}\rightarrow \mathcal{P}_{2}(\mathcal{X})$ be a curve with probability measures $p_{t}$, $v_{t}(\textbf{x}_{t})$ be the vector field which satisfy the continuity equation: $\partial_{t}p_{t} = -\texttt{div}(p_{t}v_{t})$. As shown in the following Proposition \ref{prop:decrease}, the functional $F$ is decreasing along the curve $\gamma_{t}$. Let $p^{*}$ be the optimal solution of $F$ and assume $F$ is geodesically convex, the limit of $p_{t}$ is $p^{*}$ as $t\rightarrow \infty$.

\begin{prop}
\label{prop:decrease}
$F(p_{t})$ is decreasing in time and satisfies:
\begin{align}
    \frac{\mathrm{d}F(p_{t})}{\mathrm{d}t}=-\mathbb{E}_{\textbf{x}\sim p_{t}} \lVert \nabla^{2}\varphi(\textbf{x}) ^{-1}\nabla f^{*}_{t}(\textbf{x})\rVert^{2}_{2} \leq 0
\end{align}
\end{prop}
\begin{proof}
The proof can be obtained by the following calculation:
\begin{align*}
\frac{\mathrm{d}F(p_{t})}{\mathrm{d}t}&=\int \frac{\partial F(p_{t})}{\partial p_{t}}\frac{\partial p_{t}}{\partial t}=-\int \frac{\partial F(p_{t})}{\partial p_{t}}\texttt{div}(p_{t}v_{t})\\
&= \int \langle \nabla \frac{\partial F(p_{t})}{\partial p_{t}}(\textbf{x}), v_{t}(\textbf{x})\rangle p_{t}(\textbf{x})\mathrm{d}\textbf{x}=\int \langle \nabla f^{*}_{t}(\textbf{x}), v_{t}(\textbf{x})\rangle p_{t}(\textbf{x})\mathrm{d}\textbf{x}\\
&= -\int \langle \nabla f^{*}_{t}(\textbf{x}), \nabla^{2}\varphi(\textbf{x})^{-2}\nabla f^{*}_{t}(\textbf{x})\rangle p_{t}(\textbf{x})\mathrm{d}\textbf{x}=-\int \langle \nabla f^{*}_{t}(\textbf{x}), \nabla^{2}\varphi(\textbf{x})^{-2}\nabla f^{*}_{t}(\textbf{x})\rangle p_{t}(\textbf{x})\mathrm{d} \textbf{x}\\
&=-\mathbb{E}_{\textbf{x}\sim p_{t}} \lVert \nabla^{2}\varphi(\textbf{x})^{-1} \nabla f^{*}_{t}(\textbf{x})\rVert^{2}_{2}
\end{align*}
where the equality in the first line is obtained by the definition of continuity equation, the equality in the second line is obtained by Equation (\ref{eqn:div}) and Proposition 3.1 in \cite{liu2021infinite}, and the equality in the third line is obtained by Equation (\ref{eqn:dxdt}). The proof is completed.
\end{proof}

\subsection{Relation to Stein Variational Mirror Descent \cite{shi2021sampling}}

Stein Variational Mirror Descent (\textbf{SVMD}, \cite{shi2021sampling}), most related to Algorithm \ref{alg:mirrorVT}, minimizes the KL divergence to the target distribution defined over a constrained domain by pushing particles in the dual space induced by a mirror map. In what follows, we verify that \textbf{mirrorVT} when using KL-divergence as the objective functional relates  to \textbf{SVMD} algorithms via the integral operator. First we give the definition of the integral operator as follows:
\begin{defn}
\label{def:integraloperator}
(Integral operator, \cite{rosasco2009note}) Let $k:\mathcal{X}\times\mathcal{X}\rightarrow \mathbb{R}$ be a reproducing kernel, $\mathcal{X}$ is endowed with a probability measure $p$, $L^{2}(\mathcal{X}, p)$ be the space of square integrable functions with norm $\lVert f\rVert_{p}^{2}=\langle f,f \rangle_{p}=\int |f(\textbf{x})|^{2}\mathrm{d}p(\textbf{x})$. We define $\mathcal{L}_{k,p}:L^{2}(\mathcal{X},p)\rightarrow L^{2}(\mathcal{X},p)$ to be the corresponding integral operator given by:
\begin{align}
    \mathcal{L}_{k,p}f(\textbf{x})=\int_{\mathcal{X}}k(\textbf{x},\textbf{x}^\prime)f(\textbf{x}^\prime) p(\textbf{x}^\prime) \mathrm{d}\textbf{x}
\end{align}
\end{defn}

The following theorem shows the connection between \textbf{mirrorVP} and \textbf{SVMD}.
\begin{thm}
\label{thm:relatedsvmd} Let $F(p)=\text{KL}(p||p^{*})$, for $p\in \mathcal{P}_{2}(\mathcal{X})$ be the functional to be optimized, $k$ be the reproducing kernel defined on $\mathcal{X}$, $k_{\varphi}(\textbf{y},\textbf{y}^\prime)=k(\nabla \varphi^{*}(\textbf{y}),\nabla \varphi^{*}(\textbf{y}^\prime))$ for $\textbf{y}, \textbf{y}^\prime \in \mathcal{Y}$. In $t$-th iteration, the integral operator $\mathcal{L}_{k,p_{t}}$ applied on the update direction $v_{t}$ is given by:
\begin{equation}
    \label{eqn:integraloperator}
    \mathcal{L}_{k, p_{t}} v_{t}(\textbf{x})= -\mathbb{E}_{\textbf{y}^\prime\sim q_{t}}k_{\varphi}(\textbf{y}, \textbf{y}^\prime)\nabla \log q^{*}(\textbf{y}^\prime) + \nabla_{\textbf{y}}k_{\varphi}(\textbf{y}, \textbf{y}^\prime) \text{,  }\forall \textbf{x}^\prime \in \mathcal{X}, \textbf{y}^\prime = \nabla \varphi(\textbf{x}^\prime)\in \mathcal{Y}
\end{equation}
\end{thm}

The proof of Theorem \ref{thm:relatedsvmd} is given in Appendix \ref{Appendix:relationsvmd}. We could observe that the right-hand side of (\ref{eqn:integraloperator}) is the same as that of (11) in \cite{shi2021sampling}, suggesting that in case that we optimize the KL-divergence and approximate the optimal update direction in the dual space by kernel $k$ through the integral operator, \textbf{mirrorVT} reduces to \textbf{SVMD}. Especially, when $k(\textbf{x},\textbf{x}^\prime)=1$ if $\textbf{x}=\textbf{x}^\prime$ and 0 otherwise, \textbf{mirrorVT} exactly recovers \textbf{SVMD}, so we could say that \textbf{SVMD} is a special case of \textbf{mirrorVT}. However, it might be challenging to derive a similar kernel-based update for divergences beyond KL-divergence.

\begin{algorithm}[H]
\label{alg:vfm}
\SetAlgoLined
\KwIn{Functional $F$, number of particles $N$, number of iterations: $K$, and initial weights $\textbf{w}(0)$}
\KwOut{$f_{\hat{\textbf{w}}}$}
 Stepsize $\eta=N^{-1/2}$\\
$s\leftarrow 0$\\
\While{$s< N$}{
    $\textbf{x} = \textbf{x}_{s+1}$ \tcp*[l]{pick up a sample}
    $g_{\textbf{w}(s)}(\textbf{x}) \leftarrow\nabla_{\theta}F^{*}(f_{\textbf{w}(s)}(\textbf{x}))-\nabla_{\textbf{w}}f_{\textbf{w}(s)}(\textbf{x})$ \tcp*[l]{compute gradient w.r.t network parameters}
    $\textbf{w}(s+1/2)\leftarrow \textbf{w}(s)-\eta g_{\textbf{w}(s)}(\textbf{x})$\tcp*[l]{update network parameters}
    $\textbf{w}(s+1)\leftarrow \argmin_{\textbf{w} \in \mathcal{B}^{0}}\left\{ \lVert \textbf{w} - \textbf{w}(s+1/2)\rVert_{2}\right\}$\tcp*[l]{guarantee to lie in the ball $\mathcal{B}^{0}(r_{f})$}
    $s\leftarrow s+1$
    }
$\hat{\textbf{w}} = 1/N\cdot\sum_{s=0}^{N-1}\textbf{w}(s)$
 \caption{Variational Form Maximization (\textbf{VFM}($\left\{\textbf{x}_{i} \right\}_{i\in [N]}$, $F$))}
\end{algorithm}

\subsection{Convergence Analysis}
In this section, we analyze the convergence of \textbf{mirrorVT}. Its pseudocode is shown in Algorithm \ref{alg:mirrorVT}.

In practice, we need to solve the empirical variational maximization problem (\ref{eqn:varformmax}) with finite samples via stochastic gradient descent to get the estimate $\Tilde{f^{*}_{t}}$ of the true $f^{*}_{t}$ (see Algorithm \ref{alg:vfm}).
The true Wasserstein gradient at $p_{t}$ is given by $\texttt{grad}F(p_{t})=-\texttt{div}(p_{t}\nabla f^{*}_{t})$ and its estimate is given by $-\texttt{div}(p_{t}\nabla \Tilde{f_{t}^{*}})$, so the difference is given by $\delta_{t}=-\texttt{div}(p_{t}(\nabla \Tilde{f_{t}^{*}}-\nabla f^{*}_{t}))$. Hence, the expected gradient error is defined as:
\begin{align}
    \label{eqn:graderror}
    \epsilon_{t}=\mathbb{E}\langle \delta_{t}, \delta_{t} \rangle_{p_{t}}=\mathbb{E}\int \lVert \nabla \Tilde{f_{t}^{*}}(\textbf{x})-\nabla f^{*}_{t}(\textbf{x}) \rVert^{2}_{2}p_{t}(\textbf{x})\mathrm{d}\textbf{x}
\end{align}
where the expectation is taken over the initial particles. To derive the upper bound of $\epsilon_{t}$, according to \cite{liu2021infinite}, we assume to learn the function $\Tilde{f_{t}^{*}}$ from the class $\Tilde{\mathcal{F}}$ of neural networks (\ref{eqn:nnclass}). The gradient error, defined in (\ref{eqn:graderror}), is upper bounded by: $O\left(r_{f}^{2}/N^{1/2} + r_{f}^{3}/n_{w}^{1/2} + r_{f}^{4}/n_{w}\right)$. See \cite{liu2021infinite} for the assumptions and proofs. It is noted that the gradient error decays to zero at the rate of $1/\sqrt{N}$ with a sufficiently large width $n_{w}$ of the neural network. Furthermore, the order of gradient error is independent of the iteration $t$.

We further make the following assumptions on the functionals $F$ and $G$ (see \cite{liu2021infinite} for the reference) to analyze the convergence of Algorithm \ref{alg:mirrorVT}.

\begin{algorithm}[H]
\label{alg:mirrorVT}
\SetAlgoLined
\KwIn{Functional $F:\mathcal{P}_{2}(\mathcal{X})\rightarrow \mathbb{R}$, initial measure $p_{1}\in \mathcal{P}_{2}(\mathcal{X})$, number of particles $N$, number of iterations: $T$ and stepsizes $\left\{ \eta_{t}\right\}_{t=1}^{T}$.}
\KwOut{The final set of particles $\textbf{X}_{T}$}
 Generate $N$ particles $\textbf{X}_{1}=\left\{ \textbf{x}_{1,i} \right\}_{i\in [N]}$ by drawing $N$ i.i.d. observations from $p_{1}$.\\
 \While{$t\leq T$}{
    $\Tilde{f_{t}^{*}}\leftarrow \textbf{VFM}(\textbf{X}_{t}, F)$\tcp*[l]{estimate solution of variational maximization problem}
    for $i\in [N]\text{: }$$\textbf{y}_{t,i}\leftarrow\nabla \varphi(\textbf{x}_{t,i})$ \tcp*[l]{map particles into the dual space}
    for $i\in [N]\text{: }$$v_{t}(\textbf{x}_{t,i})\leftarrow  \nabla^{2}\varphi(\textbf{x}_{t,i})^{-1}\nabla \Tilde{f_{t}^{*}}(\textbf{x}_{t,i})$\tcp*[l]{calculate update direction for each particle}
    for $i\in [N]\text{: }$$\textbf{y}_{t+1,i}\leftarrow \textbf{y}_{t,i}-\eta_{t} v_{t}(\textbf{x}_{t,i})$\tcp*[l]{push each particle in dual space}
    for $i\in [N]\text{: }$$\textbf{x}_{t+1,i} \leftarrow \nabla \varphi^{*}(\textbf{y}_{t+1,i})$\tcp*[l]{map updated particle to primal space}
    $\textbf{X}_{t+1}\leftarrow \left\{ \textbf{x}_{t+1,i}\right\}_{i\in [N]}$\tcp*[l]{update particle set}
    $t\leftarrow t+1$
    }
    \caption{The Mirrored Variational Transport (\textbf{mirrorVT}) algorithm}
\end{algorithm}

\textbf{Assumption 1.1} (Geodesic smoothness). We assume that $F$ and $G$ are geodesically $L_{1}$- and $L_{2}$-smooth with respect to the $2$-Wasserstein distance in the sense that: for $\forall p$, $p^\prime \in \mathcal{P}_{2}(\mathcal{X})$, $\forall q$, $q^\prime \in \mathcal{P}_{2}(\mathcal{Y})$:
\begin{align}
     F(p^\prime) &\leq F(p) + \langle \texttt{grad}F(p), \texttt{Exp}_{p}^{-1}(p^\prime)\rangle_{p} + \frac{L_{1}}{2}\cdot \mathcal{W}^{2}_{2}(p^\prime,p)\\
     G(q^\prime) &\leq G(q) + \langle \texttt{grad}G(q), \texttt{Exp}_{q}^{-1}(q^\prime)\rangle_{q} + \frac{L_{2}}{2}\cdot \mathcal{W}^{2}_{2}(q^\prime,q)
\end{align}
We also assume that $F$ is geodesically $\mu$-strongly convex with respect to the $2$-Wasserstein distance. That is:

\textbf{Assumption 1.2} (Geodesic strong convexity) for $\forall p$, $p^\prime \in \mathcal{P}_{2}(\mathcal{X})$, we have: 
\begin{align}
    \label{assumption:strongconvexity}
     F(p^\prime) &\geq F(p) + \langle \texttt{grad}F(p), \texttt{Exp}_{p}^{-1}(p^\prime)\rangle_{p} + \frac{\mu}{2}\cdot \mathcal{W}^{2}_{2}(p^\prime,p)
\end{align}

The geodesic strong convexity of $F$ ensures that it admits the gradient dominance defined as follows:

\textbf{Assumption 1.3} ($\mu$-gradient dominance) for $F$ satisfying (\ref{assumption:strongconvexity}) is gradient dominated in the sense that: 
\begin{align}
    \label{assumption:gradientdominance}
    \mu \cdot \left( F(p_{t}) - F(p^{*}) \right) \leq \langle \texttt{grad}F(p), \texttt{grad}F(p)\rangle_{p}
\end{align}

\begin{thm}
\label{thm:convergence}
(Convergence of Algorithm \ref{alg:mirrorVT}) Let the metric tensor $0 \prec \alpha \textbf{I}\preceq\nabla^{2}\varphi(\textbf{x})\preceq\beta \textbf{I}$ for all $\textbf{x}\in \mathcal{X}$ and under the Assumptions 1.1 and 1.2, for $t$-th iteration of Algorithm \ref{alg:mirrorVT} with the constant stepsize $\eta\in\left(0, \min\left\{ 1/4L_{1}, \alpha/h,1/\mu\beta^{2} \right\}\right)$, it holds that:
\begin{align}
\label{ineqn:convergence}
    \mathbb{E}\left[F(p_{t+1})\right] - \inf_{p\in \mathcal{P}_{2}(\mathcal{X})} F(p) \leq \rho^{t}\left( \mathbb{E}\left[F(p_{1})\right]-\inf_{p\in \mathcal{P}_{2}(\mathcal{X})} F(p)\right) + \frac{1-\rho^{t}}{1-\rho}\frac{\eta}{\alpha^{2}}\cdot  \texttt{Error}
\end{align}
where $\rho = 1-\mu\eta/2\beta^{2}\in \left[0,1 \right]$  is the contraction factor,  $\texttt{Error}=O\left(r_{f}^{2}/N^{1/2} + r_{f}^{3}/n_{w}^{1/2} + r_{f}^{4}/n_{w}\right)$ is the gradient error and the expectation is taken over the initial particle set.
\end{thm}
The proof of Theorem \ref{thm:convergence} is given in the Appendix \ref{Appendix:convergence}.
It is noted that in (\ref{ineqn:convergence}) there are optimization and gradient errors. The optimization error linearly decays by the factor of $\rho$, while statistical gradient error decays at the rate of $1/\sqrt{N}$ ($N$ is number of particles).

\section{Experiments}
In this section, we conduct a series of simulated experiments to assess algorithms for optimizing functionals with respect to the distributions defined over the constrained domain. The functionals take the forms of KL-divergence, JS-divergence and Wasserstein distance characterizing the distances to the the target distributions $p^{*}$.

\subsection{Case studies on synthetic data}
\begin{figure}
    \centering
    \begin{subfigure}{.45\textwidth}
        \centering
        \includegraphics[width=0.9\linewidth]{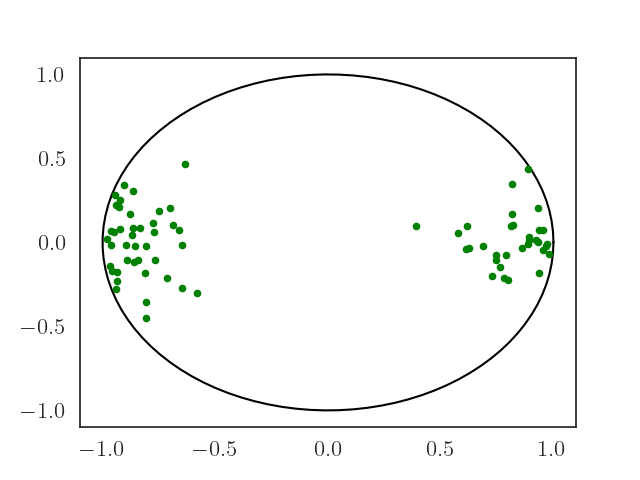}
        \caption{Samples from two-components Gaussian mixture truncated by the unit ball.}
        \label{subfig:simplex}
    \end{subfigure}
    \begin{subfigure}{.45\textwidth}
        \centering
        \includegraphics[width=0.9\linewidth]{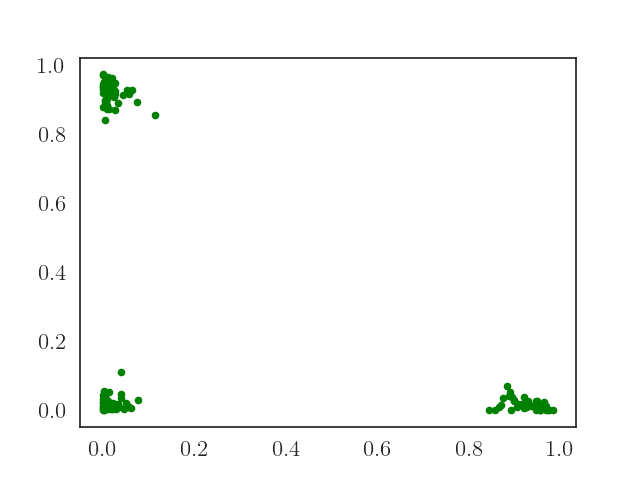}
        \caption{Samples from three-components Dirichlet mixture.}
        \label{subfig:ball}
    \end{subfigure}
    \caption{Ground truth samples.}
     \label{fig:gt}
\end{figure}

We consider two case studies with a truncated distribution and a composition distribution. For the first case, we consider an example of two-components Gaussian mixture truncated in a unit ball $\mathcal{B}^{d}=\left\{ \textbf{x}\in \mathbb{R}^{d}: \lVert \textbf{x}\rVert_{2}\leq 1 \right\}$. We generate 100 two-dimensional samples for each Gaussian component and throw away samples outside the unit ball and use the remaining ones to represent the target distribution $p^{*}$ (see Fig. \ref{fig:gt} left). The two Gaussian components are characterized by the following means and variances: $\mu_{1}=\left[-1, 0\right], \mu_{2}=\left[1, 0\right], \sigma_{1}=\sigma_{2}=0.2$.

For the second case, we consider an example of composition distribution defined on a simplex $\Delta^{d}=\left\{\textbf{x}\in \mathbb{R}_{\geq 0}^{d}: \sum_{i=1}^{d}\textbf{x}_{i}=1 \right\}$. In particular, we consider three-components Dirichlet mixture (see Fig. \ref{fig:gt} right). We generate 50 five-dimensional samples for each component. The Dirichlet distribution with unnormalized density of the form $p(\textbf{x})\propto \prod_{i=1}^{d}\textbf{x}_{i}^{\alpha_{i}-1}, \forall \textbf{x}\in \Delta^{d}$, is characterized by the concentration parameters $\alpha_{i}\ge 0,i\in [d]$. The three Dirichlet components in the simulation are characterized by $\alpha^{1}=[50,1,1,1,1]$, $\alpha^{2}=[1, 50,1,1,1]$ and $\alpha^{3}=[1,1,50,1,1]$. By visualizing only two first coordinates of data point samples, we can see that samples are mainly concentrated on three corners. For truncated Gaussian mixture, we generate 100 samples as the initial particle set which represents the initial measure $p_{1}$ (see Algorithm 2). For Dirichlet mixture, we generate 50 samples from the Dirichlet distribution characterized by $\alpha=[5,5,5,5,5]$.

\begin{figure}
    \begin{subfigure}{.32\textwidth}
        \centering
        \includegraphics[width=1.0\linewidth]{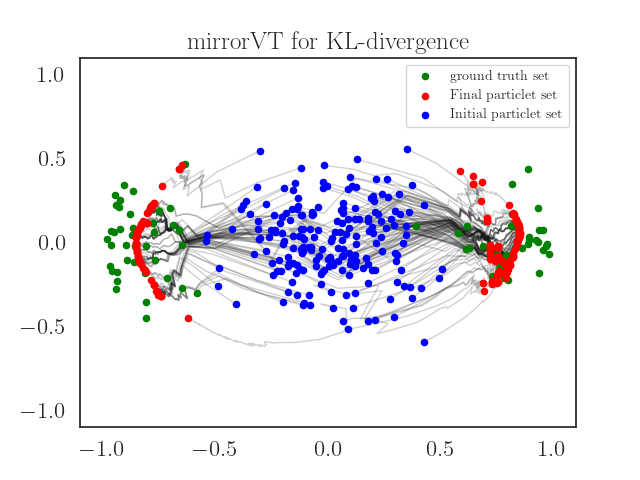}
        \label{subfig:kl}
    \end{subfigure}
    \begin{subfigure}{.32\textwidth}
        \centering
        \includegraphics[width=1.0\linewidth]{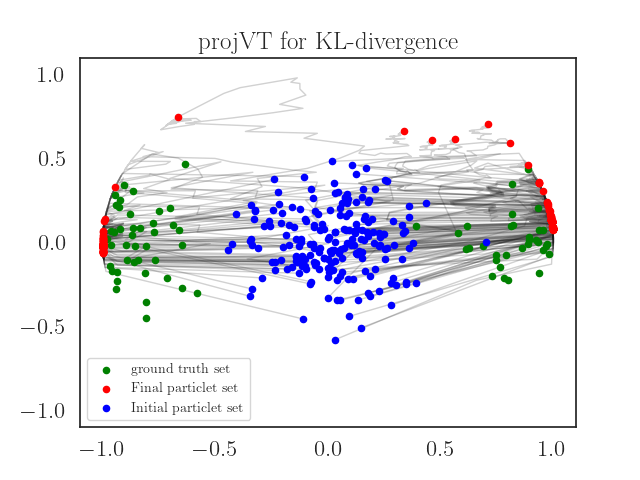}
        \label{subfig:kl}
    \end{subfigure}
    \begin{subfigure}{.32\textwidth}
        \centering
        \includegraphics[width=1.0\linewidth]{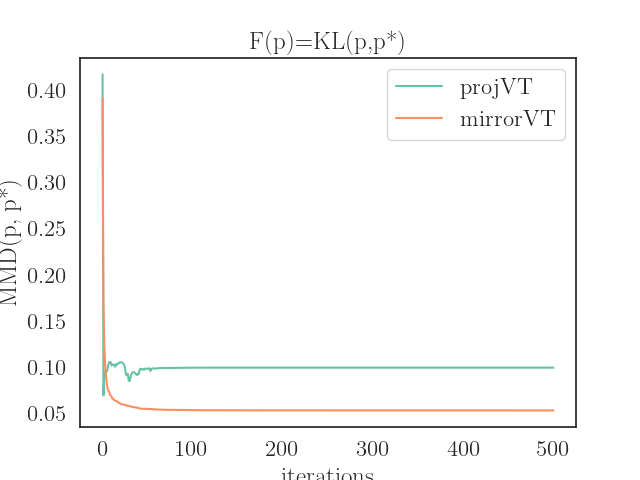}
        \label{subfig:kl}
    \end{subfigure}
    \newline
    
    \begin{subfigure}{.32\textwidth}
        \centering
        \includegraphics[width=1.0\linewidth]{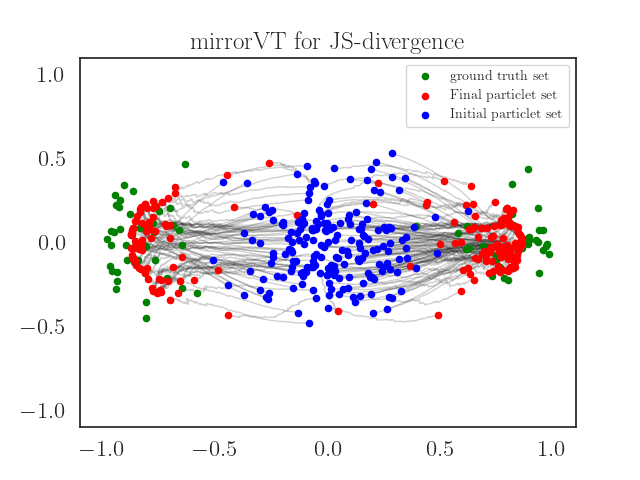}
        \label{subfig:kl}
    \end{subfigure}
    \begin{subfigure}{.32\textwidth}
        \centering
        \includegraphics[width=1.0\linewidth]{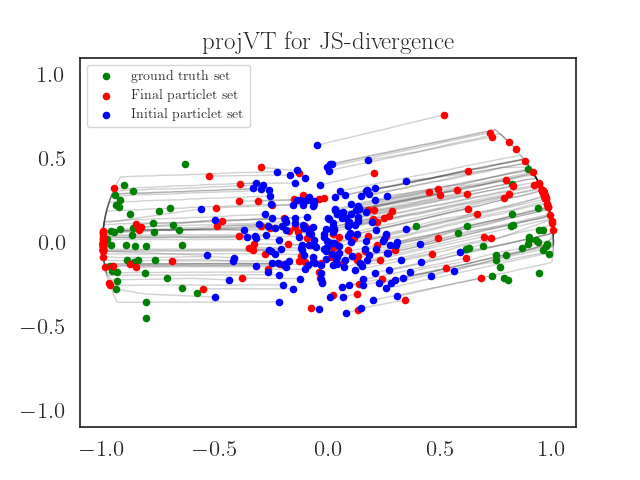}
        \label{subfig:kl}
    \end{subfigure}
    \begin{subfigure}{.32\textwidth}
        \centering
        \includegraphics[width=1.0\linewidth]{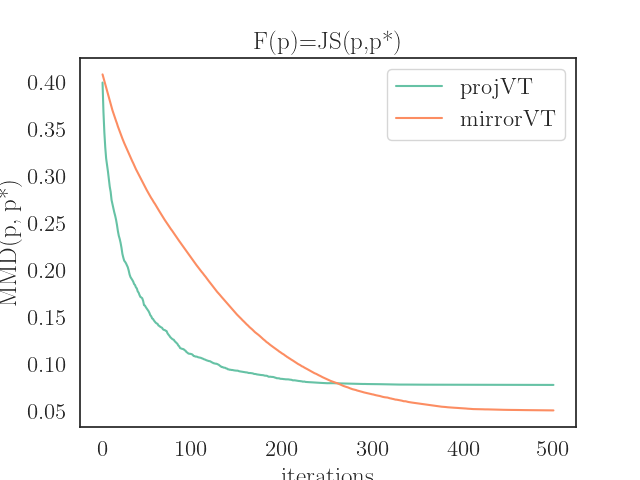}
        \label{subfig:kl}
    \end{subfigure}
    \newline
    
    \begin{subfigure}{.32\textwidth}
        \centering
        \includegraphics[width=1.0\linewidth]{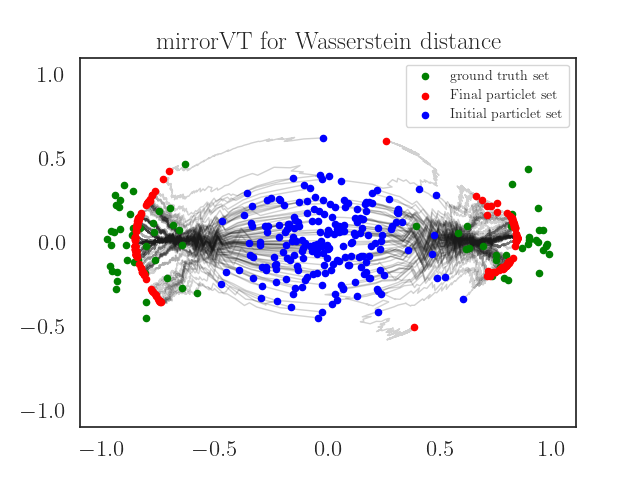}
        \label{subfig:kl}
    \end{subfigure}
    \begin{subfigure}{.32\textwidth}
        \centering
        \includegraphics[width=1.0\linewidth]{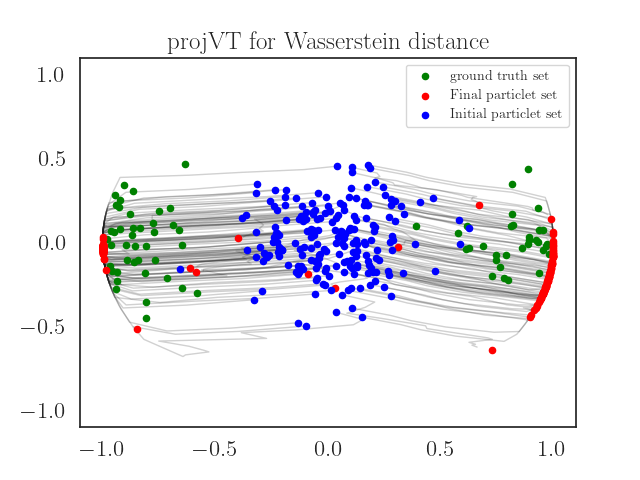}
        \label{subfig:kl}
    \end{subfigure}
    \begin{subfigure}{.32\textwidth}
        \centering
        \includegraphics[width=1.0\linewidth]{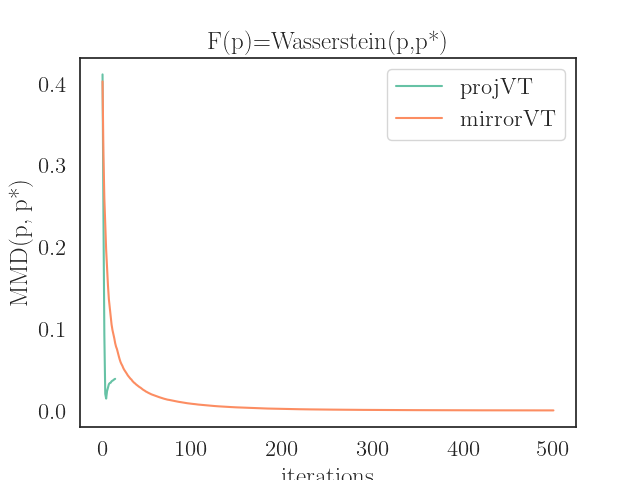}
        \label{subfig:kl}
    \end{subfigure}
    \newline
    \caption{Simulation from two-components Gaussian mixture truncated by the unit ball: trajectories of particles updated by \textbf{mirrorVT} (1st column) and \textbf{projVT} (2nd column), and the comparison of methods in terms of convergence in MMD (3rd column). $F$ takes the following forms: KL-divergence (1st row), JS-divergence (2nd row) and Wasserstein distance (3rd row). Green: particles from the target distribution, Blue: initial particle set, Red: the final updated particles.}
    \label{fig:gaussianmixture}
\end{figure}

\subsection{Comparing methods}
For comparison, we consider the \textbf{p}rojected \textbf{v}ariational \textbf{t}ransport (named \textbf{projVT}) which is a simple variant of \textbf{VT} for the constrained domain. Its update simply has two steps for $t$-th iteration: 1) pushing the particles in the direction of $\Tilde{f}^{*}_{t}$ and 2) projecting back to the constrained domain by the projection operation. In particular, for the unit ball $\mathcal{B}^{d}$, the projection operation can be easily computed by: $\texttt{proj}_{\mathcal{B}^{d}}(\textbf{x})=\argmin_{\textbf{y}\in \mathcal{B}^{d}} \lVert \textbf{x}-\textbf{y} \rVert_{2} =\textbf{x}/\max\left\{1, \lVert \textbf{x}\rVert_{2} \right\}$.

For the simplex, the projection operation can be defined as:
\begin{align*}
    \texttt{proj}_{\Delta^{d}}(\textbf{x})=\argmin_{\textbf{y}\in \Delta^{d}} \lVert \textbf{x}-\textbf{y} \rVert_{2}^{2} \text{ subject to } \sum_{i=1}^{d}\textbf{y}_{i}=1, \textbf{y}_{i}\geq 0 \forall i\in [d]
\end{align*}
which can be efficiently solved \cite{duchi2008efficient}. For \textbf{mirrorVT}, we need to define the mirror maps. In particular, for the unit ball, a possible map is:
\begin{align}
\label{eqn:unitball_mm}
    \varphi(\textbf{x})=-\log (1-\lVert \textbf{x} \rVert) - \lVert \textbf{x} \rVert
\end{align}
and for a simplex, the most natural choice of the map $\varphi$ is the entropic mirror map (see \cite{beck2003mirror}):
\begin{align}
    \label{eqn:simplex_mm}
     \varphi(\textbf{x})=\sum_{i=1}^{d-1}\textbf{x}_{i}\log \textbf{x}_{i} + \left(1 - \sum_{i=1}^{d-1}\textbf{x}_{i} \right)\log \left( 1 - \sum_{i=1}^{d-1}\textbf{x}_{i}\right), \forall \textbf{x}\in \Delta^{d}
\end{align}
The details of $\nabla \varphi(\textbf{x})$, $\nabla \varphi^{*}(\textbf{x})$ and $\nabla^{2} \varphi(\textbf{x})^{-1}$ associated with these mirror maps are given in the Appendix \ref{Appendix:mirrormap}.

\begin{figure}
    \begin{subfigure}{.32\textwidth}
        \centering
        \includegraphics[width=1.0\linewidth]{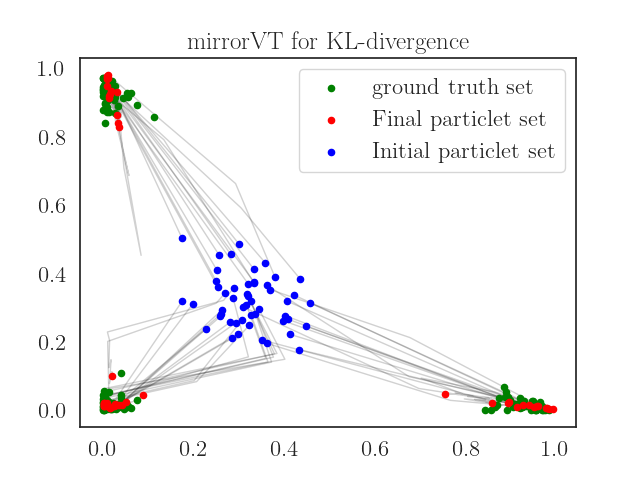}
        \label{subfig:kl}
    \end{subfigure}
    \begin{subfigure}{.32\textwidth}
        \centering
        \includegraphics[width=1.0\linewidth]{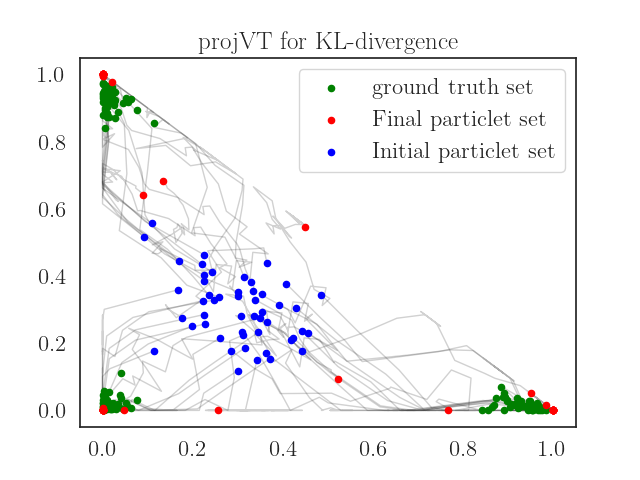}
        \label{subfig:kl}
    \end{subfigure}
    \begin{subfigure}{.32\textwidth}
        \centering
        \includegraphics[width=1.0\linewidth]{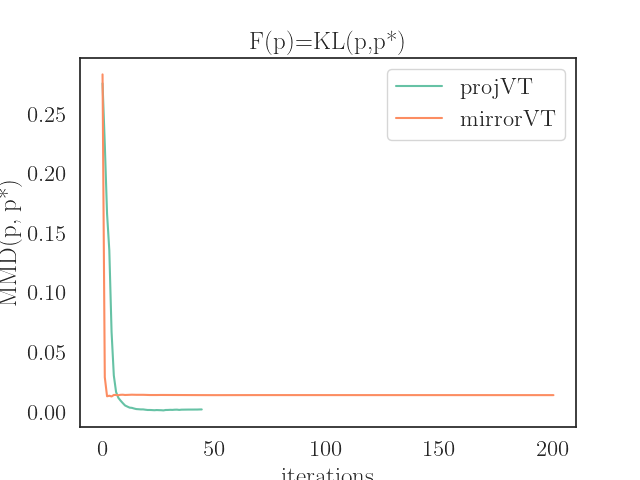}
        \label{subfig:kl}
    \end{subfigure}
    \newline

    \begin{subfigure}{.32\textwidth}
        \centering
        \includegraphics[width=1.0\linewidth]{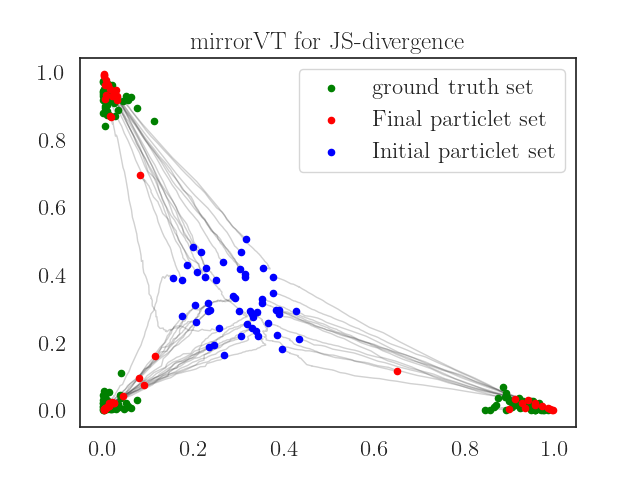}
        \label{subfig:kl}
    \end{subfigure}
    \begin{subfigure}{.32\textwidth}
        \centering
        \includegraphics[width=1.0\linewidth]{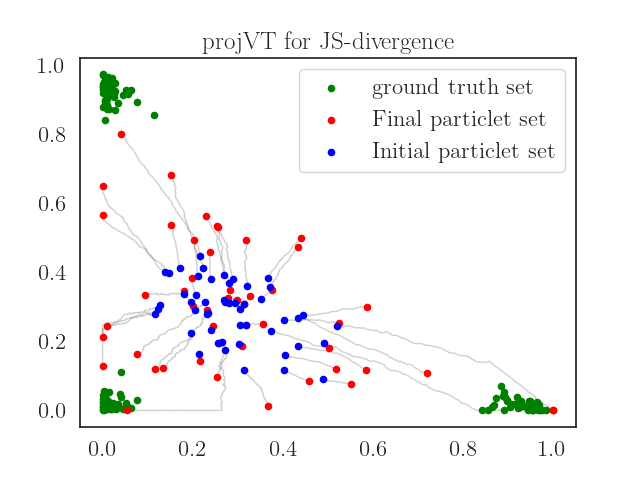}
        \label{subfig:kl}
    \end{subfigure}
    \begin{subfigure}{.32\textwidth}
        \centering
        \includegraphics[width=1.0\linewidth]{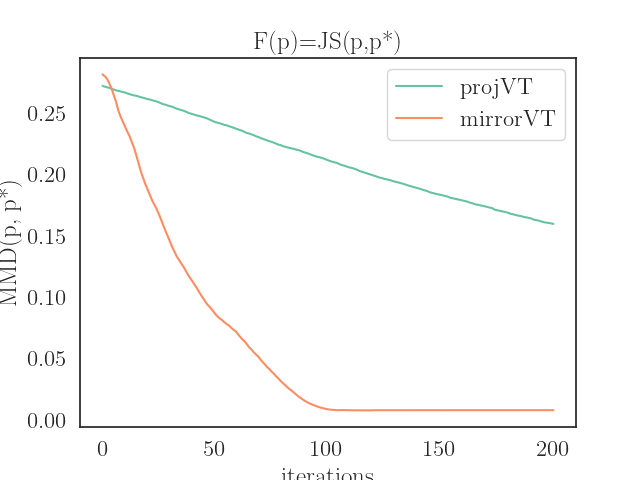}
        \label{subfig:kl}
    \end{subfigure}
    \newline
    
    \begin{subfigure}{.32\textwidth}
        \centering
        \includegraphics[width=1.0\linewidth]{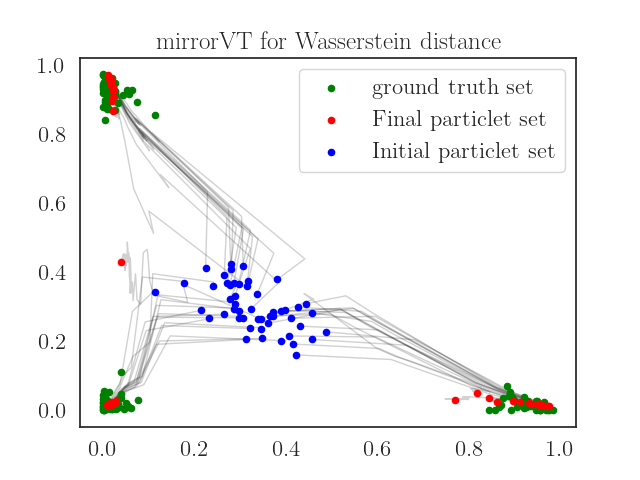}
        \label{subfig:kl}
    \end{subfigure}
    \begin{subfigure}{.32\textwidth}
        \centering
        \includegraphics[width=1.0\linewidth]{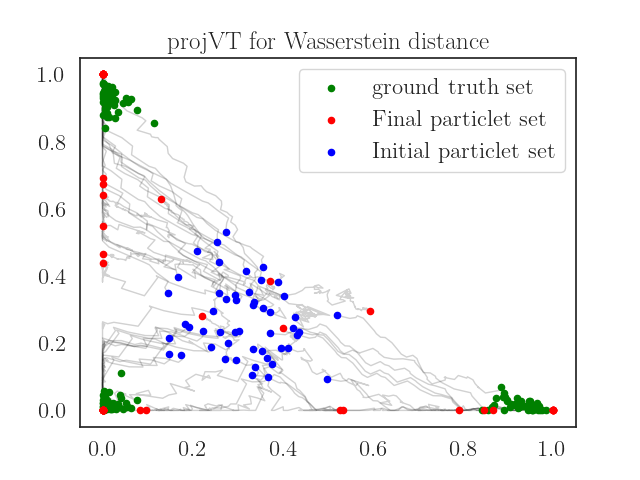}
        \label{subfig:kl}
    \end{subfigure}
    \begin{subfigure}{.32\textwidth}
        \centering
        \includegraphics[width=1.0\linewidth]{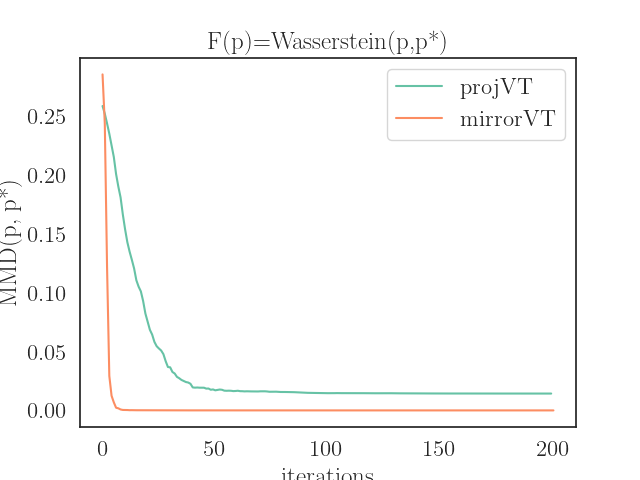}
        \label{subfig:kl}
    \end{subfigure}
    \newline
    \caption{Simulation from three-components Dirichlet mixture: trajectories of particles updated by \textbf{mirrorVT} (1st column) and \textbf{projVT} (2nd column), and the comparison of methods in terms of convergence in MMD (3rd column). $F$ takes the following forms: KL-divergence (1st row), JS-divergence (2nd row) and Wasserstein distance (3rd row). Green: particles from the target distribution, Blue: initial particle set, Red: the final updated particles.}
    \label{fig:dirichletmixture}
\end{figure}

\subsection{Results}
We compare the quality of particles updated by \textbf{projVT} and \textbf{mirrorVT} on the synthetic data samples. The step sizes for \textbf{projVT} and \textbf{mirrorVT} are fine-tuned and set to 0.01 and 0.1, respectively, and we run $T=500$ particle updates. We measure the quality of particle sets by computing their maximum mean discrepancy (MMD, \cite{gretton2012kernel}) to the ground truth particle set which represents the target distribution. We stop the particle updates if no improvement in MMD after 20 consecutive updates is observed. 
Hence we take the following forms of $F$:
KL-divergence, JS-divergence and Wasserstein distance.

On the truncated Gaussian mixture sample data, we observe that the particles updated by \textbf{mirrorVT} can well approximate the ground truth samples while those updated by \textbf{projVT} tend to stay in the boundary of the domain (see Fig. \ref{fig:gaussianmixture}). This can be explained by the nature of the projection operations. Especially, for the case of JS-divergence, \textbf{projVT} fails to converge to the target distribution. Furthermore, by measuring MMD of particles over iterations, we can find that the convergence of \textbf{mirrorVT} is significant better than that of \textbf{projVT}.

On the simplex sample data, similarly, particles updated by \textbf{mirrorVT} approximate the ground truth samples much better than those updated by \textbf{projVT} (see Fig. \ref{fig:dirichletmixture}). Again, particles by \textbf{projVT} mostly stay in the boundary with poorer quality (illustrated by the trajectories of particles and convergence of MMD). In the case of JS-divergence, \textbf{projVT} also fails to approximate the target particle set, while its particles mostly stay in the boundary in the case of Wasserstein distance.
The above observations indicate the superior performance of \textbf{mirrorVT} in the distributional optimization problem on the constrained domains.

\section{Conclusions And Future Works}
In this work, we have considered the distributional optimization problem where probability distributions are defined over the constrained domain. Motivated by the \textbf{MD} algorithm for the constrained optimization, we have presented an iterative particle-based  framework, named \textbf{mirrorVT}, to solve the problem in the unconstrained dual space induced by the mirror map. In particular, for each iteration, \textbf{mirrorVT} maps particles to the unconstrained dual space via the mirror map, then approximately performs Wasserstein gradient descent in the dual space by the \textbf{VT} framework \cite{liu2021infinite} and finally maps the updated particles back to the constrained domain. Furthermore, we also analyze theoretical properties of \textbf{mirrorVT} and characterize its convergence to the optimal solution of the objective functional. 

However, there are some open questions we would like to address in our future works, for instance, whether it is possible to accelerate the distributional optimization algorithms by applying Nesterov momentum techniques with both of the cases: unconstrained and constrained domains. We also would like to consider the development of distributional optimization algorithms for the functionals which do not admit the variational form.

\bibliographystyle{plainnat}
\bibliography{references}  






\appendix
\section{Appendices}
\subsection{Proof of Theorem \ref{thm:relatedsvmd}}.
\label{Appendix:relationsvmd}
\begin{proof}
We have assumed $F(p)=\text{KL}(p||p^{*})$ for $p,p^{*}\in \mathcal{P}_{2}(\mathcal{X})$, so $G(q)=\text{KL}(q||q^{*})$, for $q,q^{*}\in \mathcal{P}_{2}(\mathcal{Y})$. By the definition of the first variation of a functional, we have:
\begin{align*}
    \frac{\mathrm{d}}{\mathrm{d}\epsilon} G(q + \epsilon \chi)\bigg|_{\epsilon=0}=\int_{\mathcal{Y}}\frac{\partial G}{\partial q}(\textbf{y})\chi(\textbf{y})\mathrm{d}\textbf{y} \text{, for all } \chi \in \mathcal{P}_{2}(\mathcal{Y})
\end{align*}
We can compute the left-hand side as follows:
\begin{align*}
    \frac{\mathrm{d}}{\mathrm{d}\epsilon} G(q + \epsilon \chi)\bigg|_{\epsilon=0}&=\frac{\mathrm{d}}{\mathrm{d}\epsilon} \text{KL}(q + \epsilon \chi|| q^{*})\bigg|_{\epsilon=0}\\
    &=  \frac{\mathrm{d}}{\mathrm{d}\epsilon}\int (q + \epsilon \chi)\log\left(\frac{q+\epsilon \chi}{q^{*}} \right)\mathrm{d}\textbf{y}\bigg|_{\epsilon=0}\\
    &= \int \log \frac{q}{q^{*}}(\textbf{y})\chi (\textbf{y})\mathrm{d}\textbf{y}
\end{align*}
which indicates that $\partial G/\partial q= \log q - \log q^{*}$. For $t$-th iteration, the update direction $v_{t}$  is given by:
\begin{align}
\label{eqn:appendix1}
\begin{split}
    v_{t}(\textbf{x}) &= \nabla^{2}\varphi(\textbf{x})^{-1}\nabla f^{*}_{t}(\textbf{x})=\nabla g^{*}_{t}(\textbf{y})\\
    &= \nabla \log q_{t}(\textbf{y}) - \nabla \log q^{*}(\textbf{y})
\end{split}
\end{align}
for all $\textbf{x}\in \mathcal{X}, \textbf{y}=\nabla \varphi(\textbf{x})\in\mathcal{Y}$. By applying the integral operator $\mathcal{L}_{k, p_{t}}$ (see Definition \ref{def:integraloperator}) to $v_{t}$, we obtain:

\begin{align}
\begin{split}
    \mathcal{L}_{k, p_{t}} v_{t}(\textbf{x})&= \int_{\mathcal{X}}k(\textbf{x}, \textbf{x}^\prime)v_{t}(\textbf{x}^\prime)p_{t}(\textbf{x}^\prime) \mathrm{d}\textbf{x}^\prime\\
    &= \int_{\mathcal{Y}}k_{\varphi}(\textbf{y}, \textbf{y}^\prime)\nabla \log q_{t}(\textbf{y}^\prime)q_{t}(\textbf{y}^\prime)\mathrm{d}\textbf{y}^\prime-\int_{\mathcal{Y}}k_{\varphi}(\textbf{y}, \textbf{y}^\prime)\nabla \log q^{*}(\textbf{y}^\prime)q_{t}(\textbf{y}^\prime)\mathrm{d} \textbf{y}^\prime\\
    &= \int_{\mathcal{Y}}k_{\varphi}(\textbf{y}, \textbf{y}^\prime)\nabla q_{t}(\textbf{y}^\prime)\mathrm{d}\textbf{y}^\prime-\int_{\mathcal{Y}}k_{\varphi}(\textbf{y}, \textbf{y}^\prime)\nabla \log q^{*}(\textbf{y}^\prime)q_{t}(\textbf{y}^\prime)\mathrm{d}\textbf{y}^\prime\\
    &=  -\int_{\mathcal{Y}}\nabla k_{\varphi}(\textbf{y}, \textbf{y}^\prime) q_{t}(\textbf{y}^\prime) \mathrm{d}\textbf{y}^\prime-\int_{\mathcal{Y}}k_{\varphi}(\textbf{y}, \textbf{y}^\prime)\nabla \log q^{*}(\textbf{y}^\prime)q_{t}(\textbf{y}^\prime) \mathrm{d}\textbf{y}^\prime\\
    &= -\mathbb{E}_{\textbf{y}^\prime\sim q_{t}}k_{\varphi}(\textbf{y}, \textbf{y}^\prime)\nabla \log q^{*}(\textbf{y}^\prime)+ \nabla k_{\varphi}(\textbf{y}, \textbf{y}^\prime)
\end{split}
\end{align}
The first equality is obtained by the definition of the integral operator (see Definition \ref{def:integraloperator}), the second equality is obtained by using (\ref{eqn:appendix1}) and the forth equality is obtained by applying the integration by parts to the first term. The proof is completed.
\end{proof}

\subsection{Proof of Theorem \ref{thm:convergence}}
\label{Appendix:convergence}
\begin{proof}

We analyze the performance of one step of \textbf{mirrorVT}. Under the Assumption 1.1 ($L_{2}$-smoothness of $G$), for any $t\geq 0$, we have:

\begin{align}
\begin{split}
    G(q_{t+1}) &\leq G(q_{t}) + \langle \texttt{grad}G(q_{t}), \texttt{Exp}_{q_{t}}^{-1}(q_{t+1})\rangle_{q_{t}} +1/2 L_{2}\cdot \mathcal{W}^{2}_{2}(q_{t+1},q_{t})\\
    &= G(q_{t}) - \eta_{t} \langle \texttt{grad}G(q_{t}), \texttt{grad}G(q_{t})+\Tilde{\delta}_{t}  \rangle_{q_{t}}+1/2 L_{2}\eta_{t}^{2}\langle \texttt{grad}G(q_{t})+\Tilde{\delta}_{t} , \texttt{grad}G(q_{t})+\Tilde{\delta}_{t}  \rangle_{q_{t}}\\
    &=  G(q_{t}) - \eta_{t} \langle \texttt{grad}G(q_{t}), \texttt{grad}G(q_{t})  \rangle_{q_{t}} - \eta_{t}  \langle \texttt{grad}G(q_{t}), \Tilde{\delta}_{t}\rangle_{q_{t}}\\
    &+ 1/2 L_{2}\eta_{t}^{2}\langle \texttt{grad}G(q_{t})+\Tilde{\delta}_{t} , \texttt{grad}G(q_{t})+\Tilde{\delta}_{t}  \rangle_{q_{t}}
\end{split}
\end{align}
where $\eta_{t}\in \left( 0, \alpha/h\right]$ (see (11)) and $\Tilde{\delta}_{t}=-\texttt{div}(q_{t}(\nabla \Tilde{g_{t}^{*}}-\nabla g^{*}_{t}))$ is the difference between the true $2$-Wasserstein gradient at $q_{t}$ given by $\texttt{grad}G(q_{t})=-\texttt{div}(q_{t}\nabla g^{*}_{t})$ and its estimate given by $-\texttt{div}(q_{t}\nabla \Tilde{g_{t}^{*}})$. The corresponding expected gradient error for $G$ is defined as:
\begin{align}
   \Tilde{\epsilon}_{t}=\mathbb{E}\langle \Tilde{\delta}_{t}, \Tilde{\delta}_{t} \rangle_{q_{t}}=\mathbb{E}\int \lVert \nabla^{2}\varphi(\textbf{x})^{-1}\left(\nabla \Tilde{f_{t}^{*}}(\textbf{x})-\nabla f^{*}_{t}(\textbf{x})\right) \rVert^{2}_{2}p_{t}(\textbf{x})\mathrm{d} \textbf{x}
\end{align}

Also since $0 \prec \alpha \textbf{I}\preceq\nabla^{2}\varphi(\textbf{x})$ for all $\textbf{x}\in \mathcal{X}$, we have
\begin{align}
        \label{ineqn:epsilon}
     \Tilde{\epsilon}_{t} \leq \frac{1}{\alpha^{2}}\epsilon_{t}
\end{align}

By applying the basic inequality: $\langle\texttt{grad}G(q_{t}), \Tilde{\delta}_{t}\rangle \leq \frac{1}{2}\langle\texttt{grad}G(q_{t}), \texttt{grad}G(q_{t})\rangle + \frac{1}{2} \langle\Tilde{\delta}_{t}, \Tilde{\delta}_{t}\rangle $ and combining with (\ref{ineqn:epsilon}), we have:
\begin{align}
\begin{split}
    \label{inequality:1}
    G(q_{t+1}) &\leq G(q_{t}) - 1/2 \cdot\eta_{t} (1 - 2\eta_{t} L_{2})\cdot \langle\texttt{grad}G(q_{t}), \texttt{grad}G(q_{t})\rangle_{q_{t}} + \frac{\eta_{t} (1 + 2\eta_{t} L_{2})}{2\alpha^{2}} \epsilon_{t}\\
\end{split}
\end{align}

By the definition of the inner product on the tangent space and the assumption of $\mu$-strong convexity of $F$, we obtain the following inequality:
\begin{align}
\begin{split}
    \langle\texttt{grad}G(q_{t}), \texttt{grad}G(q_{t})\rangle_{q_{t}} &=\int_{\mathcal{Y}} \lVert \nabla g^{*}_{t}(\textbf{y})\rVert_{2}^{2}q_{t}(\textbf{y})\mathrm{d}\textbf{y}\\
    &=\int_{\mathcal{X}} \lVert\nabla^{2}\varphi(\textbf{x}) \nabla f^{*}_{t}(\textbf{x})\rVert_{2}^{2}p_{t}(\textbf{x})\mathrm{d} \textbf{x}\\
    &\geq \frac{1}{\alpha^{2}}\int_{\mathcal{X}} \lVert \nabla f^{*}_{t}(\textbf{x})\rVert_{2}^{2}p_{t}(\textbf{x})\mathrm{d}\textbf{x}=\frac{1}{\alpha^{2}} \langle\texttt{grad}F(p_{t}), \texttt{grad}F(p_{t})\rangle_{p_{t}}\\
    &\geq \frac{\mu}{\alpha^{2}} \left( F(p_{t})-F(p^{*})\right)
    \end{split}
\end{align}
where the first inequality is obtained by $\nabla^{2}\varphi(\textbf{x})\preceq\beta \textbf{I}$ for all $\textbf{x}\in \mathcal{X}$ and the second inequality is obtained by Assumption 1.3 (see (\ref{assumption:gradientdominance})). Thus combining (\ref{inequality:1}) and use the identity: $F(p_{t})=G(q_{t})$, we have:
\begin{align}
    F(p_{t+1}) - F(p_{t}) \leq \left[ 1 - \frac{\mu\eta_{t}}{2\beta^{2}} (1 - 2\eta_{t}L_{2}))\right]\left( F(p_{t})-F^{*}\right) + \frac{\eta_{t} (1 + 2\eta_{t} L_{2})}{2\alpha^{2}} \epsilon_{t}
\end{align}
By setting $\eta_{t} =\eta \leq \min\left\{\frac{1}{2L_{2}},\frac{1}{\mu\beta^{2}}\right\}$, we have: 
\begin{align}
    1 - \frac{\mu\eta_{t}}{2\beta^{2}} \left(1 - 2\eta_{t}L_{2})\right)\leq 1-\frac{\mu\eta_{t}}{2\beta^{2}}, 0 \leq 1-\frac{\mu\eta}{2\beta^{2}} \text{ and } \frac{\eta (1 + 2\eta L_{2})}{2\alpha^{2}} \leq \frac{\eta}{\alpha^{2}}
\end{align}

In the sequel, we define $\rho = 1-\frac{\mu\eta}{2\beta^{2}}\in \left[0,1 \right]$, we have:
\begin{align}
    F(p_{t+1}) - F(p_{t}) \leq \rho\left( F(p_{t})-F^{*}\right) + \frac{\eta_{t}}{\alpha^{2}} \epsilon_{t}
\end{align}
By forming a telescoping sequence and combining the upper bound of $\epsilon_{t}$ given in \cite{liu2021infinite}, we have:
\begin{align}
    F(p_{t+1}) - \inf_{p\in \mathcal{P}_{2}(\mathcal{X})} F(p) \leq \rho^{t}\left( F(p_{1})-\inf_{p\in \mathcal{P}_{2}(\mathcal{X})} F(p)\right) + \frac{1-\rho^{t}}{1-\rho}\frac{\eta}{\alpha^{2}}\cdot  \texttt{Error}
\end{align}
Finally, by taking the expectation over the initial particle set, we complete the proof.
\end{proof}

\subsection{Details Of The Mirror Maps}
\label{Appendix:mirrormap}
In this section, we describe more details of the mirror maps used in our simulated experiments.
\subsubsection{Mirror Map On The Unit Ball}
For the mirror map defined in (\ref{eqn:unitball_mm}), we can easily shown that:
\begin{equation*}
    \frac{\partial \varphi}{\partial \textbf{x}_{i}}=\frac{\textbf{x}_{i}}{1-\lVert \textbf{x} \rVert_{2}}, \frac{\partial \varphi^{*}}{\partial \textbf{y}_{i}}=\frac{\textbf{y}_{i}}{1+\lVert \textbf{y} \rVert_{2}},
    \frac{\partial^{2}\varphi}{\partial \textbf{x}_{i}\partial \textbf{x}_{j}}=\frac{\delta_{ij}}{1-\lVert \textbf{x} \rVert_{2}} + \frac{\textbf{x}_{i}\textbf{x}_{j}}{\lVert \textbf{x} \rVert_{2} \left( 1 - \lVert \textbf{x} \rVert_{2} \right)^{2}}
\end{equation*}
Hence the Hessian matrix can be written as: $\nabla^{2}\varphi(\textbf{x})=\frac{1}{1-\lVert \textbf{x} \rVert_{2}}\textbf{I}+\frac{1}{\lVert \textbf{x} \rVert_{2} \left( 1- \lVert \textbf{x} \rVert_{2}\right)^{2}}\textbf{x} \textbf{x}^\top$, where $\textbf{I}$ is the identity matrix. In order to obtain the inversion of Hessian matrix, we apply the celebrated Woodbury matrix identity and show that
\begin{equation*}
    \left(\nabla^{2}\varphi(\textbf{x})\right)^{-1}=\left(1 -  \lVert \textbf{x} \rVert_{2}\right) \left(\textbf{I}-\frac{1}{\lVert \textbf{x} \rVert_{2}}\textbf{x} \textbf{x}^\top \right)
\end{equation*}
\subsubsection{Mirror Map On The Simplex}
For the entropic mirror map (see \cite{beck2003mirror}), we can consider $\textbf{x}=\left[\textbf{x}_{1}, ...,\textbf{x}_{d-1}\right]\in \mathbb{R}^{d-1}$ by discarding the last entry $\textbf{x}_{d}=1 - \sum_{i=1}^{d-1}\textbf{x}_{i}$ and easily show that:
\begin{equation*}
    \frac{\partial \varphi}{\partial \textbf{x}_{i}}=\log \textbf{x}_{i} - \log \textbf{x}_{d}, \frac{\partial \varphi^{*}}{\partial \textbf{y}_{i}}=\frac{\exp{\left(\textbf{y}_{i}\right)}}{\sum_{j=1}^{d-1}\exp{\left(\textbf{y}_{j}\right)}},
    \frac{\partial^{2}\varphi}{\partial \textbf{x}_{i}\partial \textbf{x}_{j}}=\frac{\delta_{ij}}{\textbf{x}_{i}} + \frac{1}{\textbf{x}_{d}}, \forall i \in \left[d-1\right]
\end{equation*}
Hence the Hessian matrix can be written as: $\nabla^{2}\varphi(\textbf{x})=\texttt{diag}(1/\textbf{x}_{1},1/\textbf{x}_{2},...,1/\textbf{x}_{d-1})+1/\textbf{x}_{d} \textbf{1}\textbf{1}^\top$. By applying the Sherman–Morrison formula, we obtain the inverse Hessian matrix of the following form:
\begin{equation*}
    \left(\nabla^{2}\varphi(\textbf{x})\right)^{-1}=\texttt{diag}(\textbf{x}) - \textbf{x} \textbf{x}^\top
\end{equation*}
\end{document}